\documentclass[11pt]{amsart}
\usepackage[english]{babel}
\usepackage{egothic}
\usepackage[T1]{fontenc}

\usepackage{amsmath}
\usepackage{amsthm}
\usepackage{amssymb}
\usepackage[all]{xy}
\usepackage{mathrsfs}
\usepackage{enumerate}
\usepackage{physics}
\usepackage[notcite, final, notref]{showkeys}
\usepackage{dsfont}
\usepackage[dvips]{color}
\newtheorem{theorem}{Theorem}[section]
\newtheorem{problem}[theorem]{Problem}
\newtheorem{lemma}[theorem]{Lemma}
\newtheorem{proposition}[theorem]{Proposition}

\newtheorem{definition}[theorem]{Definition}

\theoremstyle{definition}
\newtheorem{remark}[theorem]{Remark}

\newtheorem{example}[theorem]{Example}

\usepackage{xcolor}

\title[]{A Hörmander-Fock space}

\author[D. Alpay]{Daniel Alpay}
\address{(DA) Schmid College of Science and Technology \\
Chapman University\\
One University Drive
Orange, California 92866\\
USA}
\email{alpay@chapman.edu}

\author[F. Colombo]{Fabrizio Colombo}
\address{(FC) Politecnico di
Milano\\Dipartimento di Matematica\\Via E. Bonardi, 9\\20133 Milano\\Italy}
\email{fabrizio.colombo@polimi.it}

\author[K. Diki]{Kamal Diki}
\address{(KD) Schmid College of Science and Technology \\
Chapman University\\
One University Drive
Orange, California 92866\\
USA}
\email{diki@chapman.edu}

\author[I. Sabadini]{Irene Sabadini}
\address{(IS) Politecnico di
Milano\\Dipartimento di Matematica\\Via E. Bonardi, 9\\20133 Milano\\Italy}
\email{irene.sabadini@polimi.it}

\author[D. C. Struppa]{Daniele C. Struppa}
\address{(DCS) Donald Bren Presidential Chair in Mathematics \\
Chapman University\\
One University Drive
Orange, California 92866\\
USA}
\email{struppa@chapman.edu}

\begin{document}
\maketitle
\begin{abstract}
In a recent paper we used a basic decomposition property of polyanalytic functions of order $2$ in one complex variable to characterize solutions of the classical $\overline{\partial}$-problem for given analytic and polyanalytic data. Our approach suggested the study of a special reproducing kernel Hilbert space that we call the Hörmander-Fock space that will be further investigated in this paper. The main properties of this space are encoded in a specific moment sequence denoted by $\eta=(\eta_n)_{n\geq 0}$ leading to a special entire function $\mathsf{E}(z)$ that is used to express the kernel function of the Hörmander-Fock space. We present also an example of a special function belonging to the class ML introduced recently by Alpay et al. and apply a Bochner-Minlos type theorem to this function, thus motivating further connections with the theory of stochastic processes.

\end{abstract}

\noindent AMS Classification: 30H20, 44A15, 46E22

\noindent Keywords: Hörmander-Fock space, moment sequence, $\overline{\partial}$-problem, special functions, ML functions
\tableofcontents
\section{Introduction}
\setcounter{equation}{0}

The research problem we are addressing in this paper can be considered as a natural continuation of the results developed in \cite{ACDSS2022_1} which focuses on a particular example involving a Gaussian function of the weighted Hörmander-type space defined in \cite[Theorem 3.2 ]{ACDSS2022_1}. 
Let us denote by $\overline{\partial}$ the classical Cauchy-Riemann operator and let $\Omega$ be a domain of $\mathbb{C}$. Consider a function $\phi$ satisfying
$$\Delta \phi(z):=\displaystyle \frac{\partial^2}{\partial z \partial \overline{z}} \phi(z)>0.$$ Hörmander's theorem \cite{hormander1966l2} asserts that given a datum $f$ satisfying suitable conditions we can solve any inhomogeneous $\overline{\partial}$-equation of the form
$$\overline{\partial} u(z,\overline{z})=f(z)$$ with a solution $u$ satisfying
$$\displaystyle \int |u(z)|^2e^{-\phi(z)}d\lambda(z)\leq \int \frac{|f(z)|^2}{\Delta \phi(z)} e^{-\phi(z)}d\lambda(z),$$
where $d\lambda(z)=dxdy$ denote the classical Lebesgue measure when $z=x+iy$.
 \\ \\
In \cite{ACDSS2022_1} the authors applied Hörmander's $L^2$-method, combined with the theory of polyanalytic functions (see \cite{abreu2014function, Balk1991, balk1997polyanalytic, vasilevski1999structure}), in order to characterize particular and general solutions of the $\overline{\partial}$-problem in one complex variable both for analytic and polyanalytic data. 
In particular we were naturally led to the study of the space
$$\mathcal{H}_p:=\left\lbrace{g\in H(\mathbb{C});\quad ||g||^2_{\mathcal{H}}:=\frac{1}{\pi}\int_{\mathbb{C}}\frac{|g(z)|^2}{(1+|z|^2)^2}e^{-p(z)}d\lambda(z)<+\infty}\right\rbrace,$$
where $p$ is a subharmonic function satisfying some suitable properties that will be described in Section $2$.
\\

In this paper we will study the special case in which  $p(z)=|z|^2$. The space $\mathcal{H}_p$ is then called Hörmander-Fock space. In this context we will introduce the special moment sequence $\eta=(\eta_n)_{n\geq 0}$ in Proposition \ref{prMS} by computing the norm of $z^n$. We will further prove various important properties of the sequence $\eta$ including the estimates proved in Proposition \ref{gammaestimate} and Proposition \ref{EST2}, and we will calculate the associated generating function, see Theorem \ref{GFS}. Finally, using the sequence $\eta$ we can introduce and study the special entire function $\mathsf{E}(z)$ presented in Definition \ref{defE}, and we can use these ideas to give a complete sequential characterization of the Hörmander-Fock space in terms of the moment sequence $\eta$, as shown in Proposition \ref{SequChara}. We also prove, see Theorem \ref{RKE}, that this space is a reproducing kernel Hilbert space whose kernel function can be expressed using the special entire function $\mathsf{E}(z)$.  \\
 \\

The results of this paper are organized as  follows: in Section 2 we review some basic facts on the $\overline{\partial}$-problem, we recall the classical Hörmander's $L^2$-estimate in one complex variable, and we review the definition of polyanalytic functions. In Section 3 we collect some results from \cite{ACDSS2022_1}  that are important for the sequel.  In Section 4 we treat the $\overline{\partial}$-problem for a given analytic datum, we give a complete characterization of the solutions to this problem and we study the Hörmander-Fock space described above.  Finally, in Section 5 we define and study a Bargmann-type kernel corresponding to the Hörmander-Fock space and present some special functions in class ML inspired from this construction.
\section{Preliminary results}
\setcounter{equation}{0}
In this section we collect different results which will be important for the sequel.
\subsection{Polyanalytic theory and special exponential integral functions }
We now recall some basic facts on polyanalytic functions of one complex variables  from \cite{Balk1991}. A polyanalytic function of order $n$ is a class $\mathcal C^n$ complex valued function $f:\Omega\subset \mathbb{C}\longrightarrow \mathbb{C}$ on a domain $\Omega$, which belongs to the kernel of the $n-th$ power, $n\geq 1,$ of the classical Cauchy-Riemann operator $\displaystyle \frac{\partial}{\partial \overline{z}}$, that is $$\displaystyle \frac{\partial^n}{\partial \overline{z}^n}f(z)=0, \quad \forall  z\in\Omega.$$
 The space of polyanalytic functions of order $n$ is denoted by $H_n(\Omega)$ and $H_1(\Omega)=H(\Omega)$.
\smallskip
An interesting fact regarding these functions is that any polyanalytic function of order $n$ can be decomposed in terms of $n$ analytic functions in fact we have a decomposition of the following form
\begin{equation}
f(z)=\displaystyle \sum_{k=0}^{n-1}\overline{z}^kf_k(z),
\end{equation}
in which all $f_k$ are analytic functions on $\Omega$.
Expanding each analytic component in series leads to an expression of the form

\begin{equation}\label{exp1}
f(z)=\displaystyle \sum_{k=0}^{n-1}\sum_{j=0}^{\infty}\overline{z}^kz^ja_{k,j},
\end{equation}
where $(a_{k,j})$ are complex coefficients.

For $n=1,2,...$ we recall that polyanalytic Fock spaces of order $n$ is defined as 
$$\mathcal{F}_n(\mathbb{C}):=\left\lbrace g\in H_n(\mathbb{C}), \quad \frac{1}{\pi}\int_{\mathbb{C}}|g(z)|^2e^{-|z|^2}d\lambda(z)<\infty \right\rbrace.$$
The reproducing kernel associated to the space $\mathcal{F}_n(\mathbb{C})$ is given by
\begin{equation}\label{Kn}
F_n(z,w)=e^{z\overline{w}}\displaystyle \sum_{k=0}^{n-1}\frac{(-1)^k}{k!}{n \choose k+1}|z-w|^{2k},
\end{equation}
for every $z,w\in\mathbb{C}.$ When $n=1$ we use the notation $\mathcal{F}(\mathbb{C})$ for the classical Fock space and $F(z,w)$ for its reproducing kernel function.\\
We recall also the special exponential integral function denoted by $E_n(z)$ and defined by (see \cite{abramowitz1964handbook})
\begin{equation}\label{Enz}
E_n(z):=\displaystyle \int_1^\infty\frac{e^{-zt}}{t^n}dt=\int_0^1e^{-\frac{z}{u}}u^{n-2}du; \quad Re(z)>0;\quad n=0, 1,...
\end{equation}
In particular we have $$\displaystyle E_0(z)=\frac{e^{-z}}{z};\quad Re(z)>0.$$
We present two important facts from \cite{abramowitz1964handbook} that will be crucial in the sequel. First of all, formula 5.1.14 of  \cite[page 229]{abramowitz1964handbook} states 
\begin{equation}\label{recFormula}
E_{n+1}(z)=\frac{1}{n}(e^{-z}-zE_n(z));\quad n=1,2, \cdots, \quad Re(z)>0.
\end{equation}
Furthermore we have
$$\frac{1}{x+n}< e^xE_n(x)\leq \frac{1}{x+n-1}; \quad x>0; \quad n=1,2,...$$

\begin{remark}\label{EnRem}
For $n=1$ we have
$$\displaystyle E_1(x)=\int_1^\infty \frac{e^{-u}}{u}du=-\gamma-\log(x)-\sum_{n=1}^\infty \frac{(-1)^n x^n}{n! n}, $$
where $\gamma$ denotes the Euler-Mascheroni constant. So, we obtain
$$\displaystyle E_1(1)=-\gamma-\sum_{n=1}^\infty \frac{(-1)^n}{n!n}.$$
\end{remark}

\subsection{Hormander's $L^2$ estimate, $\overline{\partial}$-problem and consequences}
In this work, we shall use a modified version of the Hörmander's result in the one complex variable setting as presented in \cite[Theorem 1]{berenstein1979new}(see also \cite[Theorem 1.4]{berenstein1993complex} for the case of several complex variables). Let $\Omega\subset \mathbb{C}$ be an open set, we denote by $W_p(\Omega)$ the space of all measurable functions $f:\Omega \longrightarrow \mathbb{C}$, for which there exists a constant $C>0$ such that


\begin{equation}
\displaystyle \int_{\Omega}|f(z)|^2e^{-Cp(z)}d\lambda(z)<+\infty,
\end{equation}
where $p$ is a subharmonic function satisfying:
\begin{enumerate}
\item $p(z)\geq 0$ and $\log(1+|z|^2)=O(p(z));$
\item there exist constants $C, D>0$ such that $|z -w|\leq 1$ implies that $p(z)\leq Cp(w)+D$.
\end{enumerate}
Throughout this paper, all subharmonic functions are assumed to satisfy these two properties. 
 Under these conditions Theorem 1 of \cite{berenstein1979new} can be stated as follows:
\begin{theorem}\label{HormanderBT}
Let $\Omega$ be an open subset of $\mathbb{C}$,  let $p$ be subharmonic function in $\Omega,$ and let $f$ be a function in $W_{p}(\Omega)$, such that
\begin{equation}
\displaystyle \int_{\Omega}|f(z)|^2e^{-p(z)}d\lambda(z)=M(f)<+\infty,
\end{equation}
with $d\lambda$ being the Lebesgue measure on $\Omega.$ Then, there exists a function $u$ satisfying
\begin{equation}
\displaystyle \int_{\Omega}\frac{|u(z)|^2}{(1+|z|^2)^2}e^{-p(z)}d\lambda(z)\leq \frac{M(f)}{2};
\end{equation}
and which is a solution of the equation
\begin{equation}
\overline{\partial}u=\frac{\partial}{\partial \overline{z}}u=f.
\end{equation}
Moreover, if $f\in C^{\infty}(\Omega)$ then $u\in C^{\infty}(\Omega)$.
\end{theorem}
Let us denote by $H(\Omega)$ the space of holomorphic functions on the open set $\Omega\subseteq\mathbb C$, and define
 $$\displaystyle A_p= A_p(\mathbb{C})=\left\lbrace{f\in H(\mathbb{C});\quad \exists A,B>0 \quad |f(z)|\leq A\exp(Bp(z))}\right\rbrace.$$
It is important to note that all the polynomials belong to the spaces $A_p$ because of the property (1) satisfied by $p$. These spaces are also stable by differentiation in view of property (2) of $p$. H\"ormander showed as well that
\begin{equation}\label{HormRes}
A_p(\mathbb{C})=H(\mathbb{C})\cap W_p(\mathbb{C}).
\end{equation}
Now we review some results from \cite{ACDSS2022_1} which will be useful in the sequel. 

\begin{definition}\label{Hspace}
Let $f\in W_p(\Omega)$ be such that $$\int_{\Omega}|f(z)|^2e^{-p(z)}d\lambda(z)=M(f)<+\infty.$$  We define the space
$$\mathcal{H}_{f,p}=\left\lbrace{g\in H(\Omega); \quad \int_{\Omega}\frac{|g(z)|^2}{(1+|z|^2)^2}e^{-p(z)}d\lambda(z)\leq 3M(f)}\right\rbrace.$$
\end{definition}

Applying the decomposition property of polyanalytic functions of order $2$ we have (see \cite{ACDSS2022_1} for details)

\begin{theorem}\label{H-poly}
Let $\Omega$ be an open set of $\mathbb{C}$,  let $p$ be subharmonic function in $\Omega$, and let $f$ be a function in $W_{p}(\Omega)$ such that $\overline{\partial}f=0$ and \begin{equation}
\displaystyle \int_{\Omega}|f(z)|^2e^{-p(z)}d\lambda(z)<+\infty.
\end{equation}
Then, there exists a polyanalytic function  $u$ of order $2$ which is a solution of the problem $$\overline{\partial}u=f,$$ and which can be expressed as follows
\begin{equation}\label{uexpress}
u(z,\overline{z})=\overline{z}f(z)+u_0(z), \quad
\end{equation}
where $u_0$ belongs to the function space $\mathcal{H}_{f,p}$.
\end{theorem}
\begin{remark}
The converse of Theorem \ref{H-poly} is proved in \cite[Proposition 3.3]{ACDSS2022_1} as well.
\end{remark}



As a consequence of the previous result we have 
\begin{proposition}
Let $p$ be a subharmonic function in $\mathbb{C}$.
If $f$ belongs to $A_p(\mathbb{C})$, and if
$$\displaystyle \int_{\Omega}|f(z)|^2e^{-p(z)}d\lambda(z)<+\infty,$$
  then any solution of the $\overline{\partial}$-problem considered in Theorem \ref{HormanderBT} can be expressed as follows $$u(z,\overline{z})=\overline{z}f(z)+u_0(z),$$
with $u_0\in \mathcal{H}_{f,p}(\mathbb{C})$.
\end{proposition}

We give an example inspired from the Fock space $\mathcal{F}(\mathbb{C})$.
\begin{example}
Let $p(z)=|z|^2$ and let $\Omega=\mathbb{C}$. For every fixed parameter $w\in \mathbb{C}$, we set
\begin{equation}\label{Kn}
F(z,w)=F_w(z)=e^{z\overline{w}}, \quad \forall z\in \mathbb{C}.
\end{equation}
Then, $u$ is a solution of the $\overline{\partial}-$problem with datum $F_w$ if and only if
$$u(z)=\overline{z}e^{z\overline{w}}+u_0(z),\quad \forall z\in\mathbb{C};$$
with $u_0\in\mathcal{H}_{F_w}(\mathbb{C}).$
This a direct consequence of Theorem \ref{H-poly} by taking $\Omega=\mathbb{C}$ and $f=F_w$ for every $w\in\mathbb{C}$.
\end{example}
\begin{remark}
In the previous example we have
$$M(F_w)=||F_w||_{\mathcal{F}(\mathbb{C})}^2=e^{|w|^2},\quad \forall w\in \mathbb{C}.$$
\end{remark}

\section{Reproducing kernel Hilbert space associated to $\overline{\partial}$-problem: Gaussian case}
\setcounter{equation}{0}
In this section, we study a special reproducing kernel Hilbert space induced by $\mathcal{H}_{f,p}$  where $p(z)=|z|^2$ and $\Omega=\mathbb{C}$. 
More precisely, we turn our attention to the study of a special subspace  of entire functions containing the functions $u_0$ in the formula \eqref{uexpress} of Theorem \ref{H-poly}:
\begin{definition}\label{Hpdefi}
Let $p$ be a subharmonic function on $\mathbb{C}$. We define
$$\mathcal{H}_p:=\left\lbrace{g\in H(\mathbb{C});\quad ||g||^2_{\mathcal{H}}:=\frac{1}{\pi}\int_{\mathbb{C}}\frac{|g(z)|^2}{(1+|z|^2)^2}e^{-p(z)}d\lambda(z)<+\infty}\right\rbrace.$$
and we define an inner product on $\mathcal{H}_p$ by \begin{equation}
\displaystyle \langle u,v \rangle_{\mathcal{H}_p}:=\frac{1}{\pi}\int_{\mathbb{C}}\frac{u(z)\overline{v(z)}}{(1+|z|^2)^2}e^{-p(z)}d\lambda(z).
\end{equation}
\end{definition}
\begin{remark}
The space $\mathcal{H}_p$ introduced in the previous definition corresponds to the space of functions considered by Hörmander in \cite[Theorem 2.5.3]{hormander1966l2}.
\end{remark}
\begin{definition}[Hörmander-Fock space] \label{HFS}
If we take $p(z)=|z|^2$, then the space $\mathcal{H}_{p}$ will be denoted simply by $\mathcal{H}$ and will be called the Hörmander-Fock space.
\end{definition}
\begin{proposition}
The injection $$\mathcal{F}(\mathbb{C})\hookrightarrow \mathcal{H},$$ is continuous, so that for every $g\in \mathcal{F}(\mathbb{C})$ we have
 \begin{equation}
||g||_{\mathcal{H}}\leq ||g||_{\mathcal{F}(\mathbb{C})}.
\end{equation}
\end{proposition}
\begin{proof}
For any $g\in\mathcal{F}(\mathbb{C})$ the following estimate holds $$||g||^2_{\mathcal{H}}:=\int_{\mathbb{C}}\frac{|g(z)|^2}{(1+|z|^2)^2}e^{-|z|^2}d\lambda(z)\leq  \int_{\mathbb{C}}|g(z)|^2e^{-|z|^2}d\lambda(z)=||g||_{\mathcal{F}(\mathbb{C})}<+\infty.$$
\end{proof}
\begin{remark}
In \cite{ACDSS2022_1} there is a counterexample showing that the classical Fock space $\mathcal{F}(\mathbb{C})$ and the Hörmander-Fock space $\mathcal{H}(\mathbb{C})$ are different.
\end{remark}

\begin{proposition}\label{prMS}
For every $n=0,1,...$ the norm of the monomials $z^n$ with respect to the Hilbert space $\mathcal{H}$ is given by the following moment sequence
\begin{equation}\label{MSeta}
\displaystyle \eta_n=||z^n||_{\mathcal{H}}^2=\int_0^{\infty}\frac{t^n}{(1+t)^2}e^{-t}dt\leq n!.
\end{equation}
Moreover, we have the orhogonality condition
\begin{equation}
\langle z^n,z^m \rangle_{\mathcal{H}}=0, \quad \forall n\neq m\in \mathbb{N}.
\end{equation}
\end{proposition}
\begin{proof}
We will use polar coordinates $z=re^{i\theta}$ with $r>0$ and $0\leq \theta \leq 2\pi$ so that $d\lambda(z)=rdr d\theta.$ Then, we have
\begin{align*}
\displaystyle ||z^n||_{\mathcal{H}}^2=& \frac{1}{\pi} \int_{\mathbb{C}}\frac{|z|^{2n}}{(1+|z|^2)^2}e^{-|z|^2}d\lambda(z)\\
=&\frac{1}{\pi}\int_{0}^{2\pi}\int_0^\infty \frac{r^{2n}}{(1+r^2)^2}re^{-r^2}drd\theta\\
=&2\int_0^\infty \frac{r^{2n}}{(1+r^2)^2}re^{-r^2}dr.
\end{align*}
Hence, using the change of variable $t=r^2$ so that $dt=2rdr$ we obtain

\begin{align*}
 \displaystyle
||z^n||_{\mathcal{H}}^2=& \int_0^{\infty}\frac{t^n}{(1+t)^2}e^{-t}dt\\
\leq& \int_{0}^{\infty}t^ne^{-t}dt\\
=& \Gamma(n+1)
\\
=& n!.
\end{align*}
Consider now the moment sequence \eqref{MSeta} for every $n=0,1,...$.
If $n\neq m$, using an argument similar to the one above we have

\begin{align*}
\displaystyle  \langle z^n,z^m \rangle_{\mathcal{H}}=&\frac{1}{\pi}\int_{\mathbb{C}}\frac{z^{n}\overline{z}^m}{(1+|z|^2)^2}e^{-|z|^2}d\lambda(z) \\
=& \left(\int_{0}^{2\pi}e^{(n-m)i\theta}d\theta\right)\left(\int_{0}^\infty\frac{r^{n+m}}{(1+r^2)^2}re^{-r^2}dr \right) \\
=& 0.
\end{align*}

We conclude that

$$\displaystyle \langle z^n,z^m \rangle_{\mathcal{H}}=\eta_n\delta_{n,m}, \quad \forall n,m\in \mathbb{N}.$$

\end{proof}
\begin{remark}
We can express the integral obtained in the previous computations as the Mellin transform $\mathcal{M}$ of a suitable function $f$. In fact, set $$\displaystyle f(t)=\frac{te^{-t}}{(1+t)^2}.$$
Then $$\displaystyle \mathcal{M}(f)(n):=\int_{0}^{\infty}t^{n-1}f(t)dt=\int_0^{\infty}\frac{t^n}{(1+t)^2}e^{-t}dt.$$
\end{remark}

As a first property of the moment sequence we can prove the following
\begin{proposition}\label{gammaestimate}
The moment sequence satisfies 
\begin{equation}
\frac{n!}{2^n\cdot 8}\leq \eta_n\leq n!
\end{equation}
for all $n=0,1,2,\ldots.$
\end{proposition}
\begin{proof}
It is obvious that $\eta_n\leq n!$ for every $n=0,1,...$ Then, to prove the estimate from below, we observe that $(1+t)^2=1+2t+t^2\leq 4e^t$ for every $t>0$. So, in particular
\begin{align*}
\eta_n=&\displaystyle \int_0^\infty \frac{t^n}{(1+t)^2}e^{-t}dt \\
\geq & \frac{1}{4} \int_0^\infty t^n e^{-2t}dt=\frac{1}{2^n \cdot 8}\int_0^\infty  u^ne^{-u}du.
\end{align*}
Hence, using the fact that $\int_0^\infty  u^ne^{-u}du=\Gamma(n+1)=n!$ we obtain
$$\eta_n\geq  \frac{n!}{2^n\cdot 8}; \quad \forall n\geq 0.$$
Finally, we get
$$\frac{n!}{2^n\cdot 8}\leq \eta_n\leq n!; \quad \forall n=0,1,...$$
\end{proof}

We now introduce an interesting function that encodes all the information concerning the Hörmander-Fock space $\mathcal{H}$. 
\begin{definition}\label{defE}
Given the moment sequence $\eta=(\eta_n)_{n\geq 0}$ as in \eqref{MSeta} we define the special function
\begin{equation}
\mathsf{E}(z):=\displaystyle \sum_{n=0}^\infty \frac{z^n}{\eta_n};\quad \forall z\in\mathbb{C}.
\end{equation}
\end{definition}

The next two propositions show that the special function $\mathsf{E}$ is an entire function of exponential type.

\begin{proposition}
The series $$\mathsf{E}(z)=\displaystyle \sum_{n=0}^\infty \frac{z^n}{\eta_n},$$  has radius of convergence
$$R=\lim_{n\rightarrow +\infty } (\eta_n)^{\frac{1}{n}}=\infty,$$
and therefore defines an entire function.
\end{proposition}
\begin{proof}
We observe that the radius of convergence is given by Hadamard's formula
$$\lim_{n\rightarrow \infty}(\eta_n)^{\frac{1}{n}}=\lim_{n\rightarrow \infty}e^{\frac{\log(\eta_n)}{n}}.$$
On the other hand, we note that by Proposition \ref{gammaestimate} we have $$\eta_n\geq \frac{n!}{2^n\cdot 8}; \quad \forall n\geq 0.$$
Thus, since $\lim_{n\rightarrow \infty}\frac{n!}{2^n}=+\infty$ we have also $\lim_{n\rightarrow \infty}\eta_n=+\infty.$ In particular, we have also

$$\frac{\log(\frac{n!}{2^n\cdot 8})}{n}\leq\frac{\log(\eta_n)}{n}\leq \frac{\log(n!)}{n}; \quad \forall n\geq1.$$
However, since $\lim_{n\rightarrow \infty}\frac{1}{n}\log(\frac{n!}{2^n\cdot 8})=+\infty$ we obtain that the radius of convergence is given by  $$R=\lim_{n\rightarrow \infty}(\eta_n)^{\frac{1}{n}}=\lim_{n\rightarrow \infty}e^{\frac{\log(\eta_n)}{n}}=+\infty.$$
\end{proof}

As to the growth of $\mathsf{E}(z)$ we prove the following result:
\begin{proposition}\label{EstimaeE}
For every $z\in \mathbb{C}$, it holds that
\begin{equation}
e^{|z|} \leq \mathsf{E}(|z|)\leq 8 e^{2|z|}.
\end{equation}
\end{proposition}
\begin{proof}
First of all, for every $z\in\mathbb{C}$ observe that $$\displaystyle e^{|z|}=\sum_{n=0}^{\infty}\frac{|z|^n}{n!}\leq \sum_{n=0}^{\infty}\frac{|z|^n}{\eta_n}=\mathsf{E}(|z|).$$
On the other hand, using Proposition \ref{gammaestimate} we have $$\frac{n!}{2^n\cdot 8}\leq\eta_n; \quad \forall n\geq 0.$$ Hence, we have $$\frac{1}{\eta_n}\leq 8 \frac{2^n}{n!}; \quad \forall n\geq 0.$$
Thus, for every $z\in \mathbb{C}$ it follows
$$\mathsf{E}(|z|)= \sum_{n=0}^{\infty}\frac{|z|^n}{\eta_n}\leq 8 \sum_{n=0}^{\infty}\frac{(2|z|)^n}{n!}=8e^{2|z|}.$$
Hence, we have

$$e^{|z|} \leq \mathsf{E}(|z|)\leq 8 e^{2|z|}; \quad \forall z\in \mathbb{C}.$$
\end{proof}

Now we investigate other interesting properties related to the moment sequence $\eta=(\eta_n)_n$ and the special function $\mathsf{E}(z)$:
\begin{proposition}
Let us consider the moment sequence given by \eqref{MSeta}.
Then 
\begin{equation}
\displaystyle \sum_{n=0}^\infty\frac{\eta_n}{n!}=1.
\end{equation}

\end{proposition}
\begin{proof}
Starting from the definition of the moment sequence $\eta_n$ we note that, for every $N=0,1,...$,
\begin{align*}
\displaystyle \sum_{n=0}^N \frac{\eta_n}{n!}&=\int_0^\infty\sum_{n=0}^N\frac{1}{n!}\frac{t^n}{(1+t)^2}e^{-t}dt\\
&=\int_0^\infty f_N(t)dt;\\
\end{align*}
where we have set $$f_N(t):=\sum_{n=0}^N\frac{1}{n!}\frac{t^n}{(1+t)^2}e^{-t};\quad N\geq 0, t>0.$$
First of all, it is clear that $$\displaystyle \lim_{N\rightarrow \infty}f_N(t)=\frac{1}{(1+t)^2}, \quad \forall t>0.$$
Moreover, it is easy to check that for every $N\geq 0$, and every $t>0$, we have
$$f_{N+1}(t)-f_{N}(t)=\frac{t^{N+1}}{(N+1)!}\frac{e^{-t}}{(1+t)^2}\geq 0.$$
So, the sequence $(f_N(t))_{N\geq 0}$ is an increasing sequence of positive measurable functions which converges to $\frac{1}{(1+t)^2}$ and thus using Beppo-Levi's Theorem we obtain

\begin{align*}
\displaystyle \sum_{n=0}^\infty \frac{\eta_n}{n!}&=\lim_{N\rightarrow \infty}\sum_{n=0}^N \frac{\eta_n}{n!} \\
&=\lim_{N\rightarrow \infty} \int_0^\infty f_N(t)dt\\
&=\int_0^\infty\lim_{N\rightarrow \infty}f_N(t)dt\\
&= \int_0^\infty \frac{1}{(1+t)^2} dt\\
&= 1.\\
\end{align*}

\end{proof}
We now define the following generating function:
\begin{theorem}\label{GFS}
For every $z\in \mathbb{C}$ such that $Re(z)>-1$, we have
\begin{equation}
\displaystyle \sum_{n=0}^\infty (-1)^n\frac{\eta_n}{n!}z^n=1-(z+1)e^{z+1}E_1(z+1).
\end{equation}
\end{theorem}
\begin{proof}
Let $z\in\mathbb{C}$ be such that $Re(z)>-1$, and let $\eta=(\eta_n)_n$ be the moment sequence given by \eqref{MSeta}. We consider the function $S(z)$ defined by the series 
\begin{equation}
\displaystyle S(z)=\sum_{n=0}^{+\infty} (-1)^n\eta_n\frac{z^n}{n!}.
\end{equation}
By the Lebesgue theorem we have
\begin{align*}
\displaystyle S(z)&=\sum_{n=0}^{+\infty} (-1)^n\frac{z^n}{n!}\int_0^\infty \frac{t^n}{(1+t)^2}e^{-t}dt \\
&=\int_0^\infty \left(\sum_{n=0}^{+\infty} \frac{(-zt)^n}{n!} \frac{e^{-t}}{(1+t)^2}\right)dt\\
&=\int_0^\infty e^{-zt}\frac{e^{-t}}{(1+t)^2}dt\\
&=\int_0^\infty \frac{e^{-(z+1)t}}{(1+t)^2}dt.\\
\end{align*}
Then, making the change of variables $u=1+t$, we obtain
\begin{align*}
\displaystyle S(z)&=\int_1^\infty e^{-(z+1)(u-1)}\frac{du}{u^2}\\
&=e^{z+1}\int_1^\infty e^{-(z+1)u}\frac{du}{u^2}.\\
\end{align*}
Now we compute the term $I(z)=\displaystyle\int_1^\infty e^{-(z+1)u}\frac{du}{u^2}$. To this end, we integrate by part setting $u_0=e^{-(z+1)u}$ and $v_0'=u^{-2}$ so that we have $u_0'=-(z+1)e^{-(z+1)u}$ and $v_0=-u^{-1}$. Thus  for $Re(z+1)>0$, we get
$$I(z)=e^{-(z+1)}-\int_1^\infty (z+1)u^{-1}e^{-(z+1)u}du.$$
We apply a second integration by parts setting $u_0=e^{-(z+1)u}, $ $v_0'=u^{-1}$ so that we have $u_0'=-(z+1)e^{-(z+1)u}$ and $v_0=\log(u)$. We obtain
\begin{align*}
I(z)&=e^{-(z+1)}-(z+1)\left(\left[\log(u)e^{-(z+1)u}\right]_1^\infty+\int_1^\infty (z+1)e^{-(z+1)u}\log(u) du \right)\\
&=e^{-(z+1)}-(z+1)\left(0+(z+1)\int_1^\infty \log(u) e^{-(z+1)u}du \right)\\
&=e^{-(z+1)}-(z+1)^2\int_1^\infty \log(u) e^{-(z+1)u}du.\\
\end{align*}
We now use the  formula (see \cite[page 571]{zwillinger2007table})
\begin{equation}
\displaystyle \int_1^\infty e^{-\mu x}\log(x)dx=-\frac{1}{\mu}Ei(-\mu), \quad  Re(\mu)>0,
\end{equation}
where $Ei(x):=-\int_{-x}^{+\infty}\frac{e^{-t}}{t}dt$ denotes the exponential integral, see \cite{abramowitz1964handbook, zwillinger2007table} for more details.
Taking $\mu=z+1$ we have $Re(\mu)>0$ and $$\int_1^\infty \log(u) e^{-(z+1)u}du=-\frac{1}{z+1}Ei(-(z+1)).$$ Moreover, using the well-known fact that $E_1(\mu)=-Ei(-\mu)$ we obtain $$I(z)=e^{-(z+1)}-(z+1)^2\frac{E_1(z+1)}{z+1},$$
and this leads to
$$S(z)=(z+1)I(z)=1-(z+1)E_1(z+1), \quad z\in\mathbb{C}, \quad Re(z)>-1,$$
which completes the proof.
\end{proof}
\begin{remark}
Note the implicit role of the Laplace transform $\mathcal{L}$ in the previous calculations. In fact, we observe that for every $z\in \mathbb{C}$ such that $Re(z)>-1$, we have
\begin{align*}
\displaystyle \sum_{n=0}^\infty (-1)^n\frac{\eta_n}{n!}z^n&=\int_0^\infty \left(\sum_{n=0}^\infty \frac{t^n}{n!}(-z)^n \right)\frac{e^{-t}}{(1+t)^2}dt\\
&=\int_0^\infty e^{-zt} \frac{e^{-t}}{(1+t)^2}dt\\
&=\int_0^\infty e^{-(z+1)t}\frac{1}{(1+t)^2}dt\\
&=\mathcal{L}\left(\frac{1}{(1+t)^2}\right)(z+1)\\
&=1-(z+1)e^{z+1}E_1(z+1).\\
\end{align*}
\end{remark}
So using the induction formula \eqref{recFormula} we can give an expression of the moment sequence $(\eta_n)_{n\geq 0}$ in terms of the special functions $(E_n)_{n\geq 0}$. To this end, we need first to prove a technical lemma:
\begin{lemma}\label{techlem}
For every $n=0,1\cdots$ we have
\begin{equation}
 \int_1^\infty \frac{(u-1)^n}{u}e^{-u}du= \Gamma(n+1)E_{n+1}(1).
\end{equation}
\end{lemma}
\begin{proof}
For $n=0$ it is clear that
$ \int_1^\infty \frac{e^{-u}}{u}du=\Gamma(1)E_1(1).$ So, we assume by induction that the result holds for a certain $n\geq 1$ and we prove it for $n+1$. Indeed,  thanks to the induction hypothesis we have
\begin{align*}
\displaystyle \int_1^\infty \frac{(u-1)^{n+1}}{u}e^{-u}du&=\int_1^\infty (u-1)^ne^{-u}du-\int_1^\infty \frac{(u-1)^n}{u}e^{-u}du\\
&=\int_1^\infty (u-1)^ne^{-u}du-\Gamma(n+1)E_{n+1}(1)\\
\end{align*}
Then, using a simple change of variables we can show 
$$\displaystyle \int_1^\infty (u-1)^{n}e^{-u}du=\int_0^\infty t^ne^{-t-1}dt=e^{-1}\Gamma(n+1).$$
However thanks to \eqref{recFormula} for $z=1$  we observe
\begin{equation}
nE_{n+1}(1)=e^{-1}-E_n(1), \quad n=1,2, \cdots
\end{equation}

Hence, we conclude
\begin{align*}
\displaystyle \int_1^\infty \frac{(u-1)^{n+1}}{u}e^{-u}du&=\Gamma(n+1)\left(e^{-1}-E_{n+1}(1)\right)\\
&=(n+1)\Gamma(n+1)E_{n+2}(1)\\
&=\Gamma(n+2)E_{n+2}(1).\\
\end{align*}
\end{proof}

\begin{proposition}\label{expression2}

We have $\eta_0=1-e E_1(1)$. Moreover,
for every $n=1,2,...$ we have
\begin{equation}\label{ga2}
\eta_n=(e(1+n)E_n(1)-1)\Gamma(n).
\end{equation}

\end{proposition}
\begin{proof}
We observe that the change of variables $u=1+t$ leads to $$\displaystyle \eta_0=\int_0^\infty \frac{e^{-t}}{(1+t)^2}dt= e\int_{1}^{\infty} \frac{e^{-u}}{u^2}du=eE_2(1)=1-eE_1(1).$$
We use induction. Assume expression \eqref{ga2} holds for a given $n$; we will prove it holds for $n+1$. In fact, it is easy to see using the change of variables $u=1+t$ that we have
$$\eta_n=e \int_1^\infty \frac{(u-1)^n}{u^2}e^{-u}du$$
We observe also that applying Lemma \ref{techlem} we have
\begin{equation}
\displaystyle \int_1^\infty \frac{(u-1)^n}{u}e^{-u}du=\Gamma(n+1)\Gamma(-n,1)=\Gamma(n+1)E_{n+1}(1),
\end{equation}
where $\Gamma(s,x)$ denotes the incomplete Gamma function.
We use the relation between $E_n$ and $E_{n+1}$ which leads to
$$E_n(1)=e^{-1}-nE_{n+1}(1).$$

 Thus, developing the calculations we have
\begin{align*}
\displaystyle \eta_{n+1}&=e \int_1^\infty \frac{(u-1)^{n+1}}{u^2}e^{-u}du\\
&=e \left( \int_1^\infty \frac{(u-1)^n}{u}e^{-u}du -\int_1^\infty \frac{(u-1)^n}{u^2}e^{-u}du\right)\\
&=e\int_1^\infty \frac{(u-1)^n}{u}e^{-u}du-\eta_n.\\
\end{align*}
At this stage we can use the induction hypothesis combined with  Lemma \ref{techlem} to obtain
\begin{equation}\label{F1}
\displaystyle \eta_{n+1}=e \Gamma(n+1)E_{n+1}(1)-(e(1+n)E_{n}(1)-1)\Gamma(n).
\end{equation}
However, since $E_n(1)=e^{-1}-nE_{n+1}(1)$ we deduce
\begin{equation}\label{F2}
e(1+n)E_n(1)-1=n\left(1-e(1+n)E_{n+1}(1)\right).
\end{equation}
Thus, we insert the expression \eqref{F2} in formula \eqref{F1} and use the fact $\Gamma(n+1)=n\Gamma(n)$ to get
\begin{align*}
\displaystyle \eta_{n+1}&=e\Gamma(n+1)E_{n+1}(1)-n\Gamma(n)\left(1-e(1+n)E_{n+1}(1) \right)\\
&=\Gamma(n+1)\left(eE_{n+1}(1)-1+e(1+n)E_{n+1}(1)\right)\\
&=\Gamma(n+1)\left((2+n)eE_{n+1}(1)-1 \right).\\
\end{align*}

Finally, we deduce that
\begin{equation}
\eta_{n+1}= \left[e(2+n)E_{n+1}(1)-1\right]\Gamma(n+1).
\end{equation}
This ends the proof.

\end{proof}
\begin{proposition}
The moment sequence $(\eta_n)_{n\geq 0}$ can be expressed as follows $$\displaystyle \eta_n=e\sum_{k=0}^{n}(-1)^{n-k}{n\choose k}E_{2-k}(1),\quad n=0,1,\cdots$$
\end{proposition}
\begin{proof}
We observe that
\begin{align*}
\displaystyle \eta_n&=e \int_1^\infty \frac{(u-1)^n}{u^2}e^{-u}du\\
&=e\sum_{k=0}^{n}(-1)^{n-k}{n \choose k}\int_1^\infty u^{k-2}e^{-u}du\\
&=e\sum_{k=0}^{n}(-1)^{n-k}{n \choose k} \Gamma(k-1,1).\\
\end{align*}
We now use the relation between the incomplete Gamma function $\Gamma(a,z)$ and the special functions $E_n(z)$ given by (see formula 5.1.45 from \cite{abramowitz1964handbook} )

$$E_n(z)=z^{n-1}\Gamma(1-n,z).$$
So, taking $z=1$ we have $\Gamma(k-1,1)=E_{2-k}(1),$ which gives
$$\displaystyle \eta_n=e\sum_{k=0}^{n}(-1)^{n-k}{n\choose k}E_{2-k}(1).$$
\end{proof}
In the next result we obtain a better estimate for the moment sequence $(\eta_n)_{n\geq 0}$. 
\begin{proposition}\label{EST2}
For every $n=1,2,...$ we have
\begin{equation}
0< \eta_n\leq \frac{\Gamma(n)}{n}.
\end{equation}
\end{proposition}
\begin{proof}
We note that by \cite[page 229, Formula 5.1.19]{abramowitz1964handbook} we have the following inequality for the integral exponential function $E_n(x)$:
$$\frac{1}{x+n}< e^xE_n(x)\leq \frac{1}{x+n-1}; \quad x>0; \quad n=1,2,...$$
So, setting $x=1$ in the previous relation we get
\begin{equation}
\frac{1}{1+n}< eE_n(1)\leq \frac{1}{n}; \quad n=1,2,...
\end{equation}
Hence, it follows that
$$0< e(1+n)E_n(1)-1\leq \frac{1}{n}; \quad n=1,2,...$$

We multiply by $\Gamma(n)$ the previous inequality and use the expression of the moment sequence $\eta_n$ obtained in Proposition \ref{expression2} to get

$$0<\eta_n \leq \frac{\Gamma(n)}{n}; \quad n=1,2,...$$
\end{proof}


\begin{remark}
We note that Remark \ref{EnRem} yields that
$$\displaystyle \eta_1=-2e \gamma-2e \sum_{k=1}^\infty \frac{(-1)^k}{k\cdot k!}-1.$$
It is possible to show also that for every $a>0$ we have
$$\displaystyle \eta_1=-1+2e a\int_0^\infty \frac{e^{-at}}{1+at}dt.$$
These calculations are based on properties of the special functions $E_n$, see e.g. the book \cite[page 230]{abramowitz1964handbook}.
\end{remark}
\begin{proposition}
There exists a constant $C>0$ such that for every $z\in \mathbb{C}$ and $f\in \mathcal{H}$ we have
\begin{equation}
|f(z)|\leq \sqrt{\mathsf{E}(|z|^2)}||f||_{\mathcal{H}}\leq Ce^{|z|^2}||f||_{\mathcal{H}}.
\end{equation}
\end{proposition}
\begin{proof}
We use a standard argument based on the Cauchy-Schwarz inequality. The second estimate is based on Proposition \ref{EstimaeE}.
\end{proof}
\begin{proposition}\label{ONB}
The family of functions defined by
\begin{equation}
e_n(z)=\displaystyle \frac{z^n}{\sqrt{\eta_n}}, \quad n\in\mathbb{N}, z\in \mathbb{C},
\end{equation}
forms an orthonormal basis of the Hörmander-Fock space $\mathcal{H}$.
\end{proposition}
\begin{proof}
Let $f\in\mathcal{H}$. Since $f$ is entire it is clear that we can write $$\displaystyle f(z)=\sum_{n=0}^{\infty}e_n(z)\alpha_n,$$
with complex coefficients $(\alpha_n)_n$. We use the computations developed in Proposition \ref{prMS} to claim that
$$\langle e_n, e_m \rangle_{\mathcal{H}}=\delta_{n,m},\quad n,m=0,1,\cdots $$
Moreover, if $$\langle f,e_n\rangle_{\mathcal{H}}=0, \quad \forall n=0,1,\cdots$$ we obtain that all the coefficients $\alpha_n$ vanish, which implies that $f=0$ everywhere. Hence $(e_n)_{n\geq 0}$ form an orthonormal basis of $\mathcal{H}$.
\end{proof}

\begin{proposition}\label{SequChara}
An entire function $g(z)=\displaystyle \sum_{n=0}^\infty z^na_n$ belongs to $\mathcal{H}$ if and only if
\begin{equation}
\displaystyle \sum_{n=0}^\infty \eta_n |a_n|^2<+\infty.
\end{equation}
\end{proposition}
\begin{proof}

We apply the orthogonality conditions combined with the following fact
\begin{align*}
\displaystyle ||g||^2_\mathcal{H}&= \frac{1}{\pi}\int_{\mathbb{C}}\frac{|g(z)|^2}{(1+|z|^2)^2}e^{-|z|^2}d\lambda(z)\\
&=\sum_{n,m=0}^{\infty} a_n\overline{a_m}\langle z^n, z^m \rangle_{\mathcal{H}}\\
&=\sum_{n=0}^{\infty} \eta_n |a_n|^2<\infty.\\
\end{align*}

\end{proof}

\begin{theorem}\label{RKE}
The Hörmander-Fock space $\mathcal{H}$ is a reproducing kernel Hilbert space with a reproducing kernel given by
\begin{equation}\label{Ker}
\mathsf{K}(z,w):=\displaystyle \sum_{n=0}^{\infty}\frac{z^n\overline{w}^n}{\eta_n}=\mathsf{E}(z\overline{w}),\quad \forall z,w\in \mathbb{C}.
\end{equation}
Moreover, for every fixed $z\in\mathbb{C}$ we have
\begin{equation}
||\mathsf{K}_z||_\mathcal{H}^2=\mathsf{K}(z,z)=\mathsf{E}(|z|^2).
\end{equation}
\end{theorem}
\begin{proof}
Since $\lbrace e_n(z);\quad n=0,1,\cdots \rbrace$ forms an orthonormal basis of $\mathcal{H}$, the reproducing kernel of $\mathcal{H}$ can be computed  for every $z,w \in \mathbb{C}$ as follows

\begin{align*}
\displaystyle \mathsf{K}(z,w):&=\sum_{n=0}^{\infty}e_n(z)\overline{e_n(w)} \\
&=\sum_{n=0}^{\infty}\frac{z^n\overline{w}^n}{\eta_n}\\
&=\mathsf{E}(z\overline{w}).\\
\end{align*}
In particular, it is clear that
$$\mathsf{K}(z,z)=\mathsf{E}(|z|^2), \quad \forall z\in \mathbb{C}.$$

\end{proof}
\section{Bargmann-type kernel and ML class examples }
\setcounter{equation}{0}

We now consider a Segal-Bargmann type kernel associated to the moment sequence $(\eta_n)_{n \geq 0}$, which can be defined using the normalized Hermite functions $(\psi_n)_{n\geq 0}$ as follows:
\begin{definition}
 We define the $\eta-$Bargmann kernel using the generating function below
\begin{equation}
\mathsf{A_\eta}(z,x)=\displaystyle \sum_{n=0}^{\infty}\frac{z^n}{\sqrt{\eta_n}}\psi_n(x);\quad z\in\mathbb{C}, x\in\mathbb{R}.
\end{equation}
\end{definition}
\begin{proposition}
For every $z\in \mathbb{C}$, the function $A_z: x\longrightarrow A_\eta(z,x)$ belongs to the Hilbert space $L^2(\mathbb{R})$. Moreover we have
$$||A_z||_{L^2(\mathbb{R})}^2=\mathsf{E}(|z|^2).$$
\end{proposition}
\begin{proof}
We omit the proof since it is based on standard calculations using the orthogonality of the Hermite functions and the definition of the function $\mathsf{E}(|z|^2)$.
\end{proof}
It is possible to consider another generating function using the normalized Hermite functions $(\psi_n(x))_{n\geq 0}$:
\begin{proposition}
For every $z\in\mathbb{C}$ and $x\in \mathbb{R}$ we have
\begin{equation}
\displaystyle \sum_{n=0}^\infty \frac{\eta_n}{\sqrt{n!}}z^n\psi_n(x)=e^{-\frac{x^2}{2}}\int_0^\infty\frac{e^{-\frac{z^2t^2}{2}+(\sqrt{2}zx-1)t}}{(1+t)^2}dt
\end{equation}
\end{proposition}
\begin{proof}
The computations are based on the  definition of the moment sequence $(\eta_n)_{n\geq 0}$ and the use of the classical generating function associated to normalized Hermite functions  given by
$$\displaystyle \sum_{n=0}^\infty \frac{z^n}{\sqrt{n!}}\psi_n(x)=e^{-\frac{1}{2}(z^2+x^2)+\sqrt{2}zx}, \quad \forall z\in\mathbb{C},\forall x\in \mathbb{R}.$$
In fact, we have
\begin{align*}
\displaystyle  \sum_{n=0}^\infty \frac{\eta_n}{\sqrt{n!}}z^n\psi_n(x)&=\sum_{n=0}^{\infty} \frac{z^n}{\sqrt{n!}}\left(\int_0^\infty \frac{t^n}{(1+t)^2}e^{-t}dt\right)\psi_n(x)\\
&=\int_0^\infty \left( \sum_{n=0}^{\infty}\frac{(tz)^n}{\sqrt{n!}}\psi_n(x)\right)\frac{e^{-t}}{(1+t)^2}dt\\
&=\int_0^\infty\frac{e^{-\frac{1}{2}(z^2t^2+x^2)+(\sqrt{2}zx-1)t}}{(1+t)^2}dt.\\
\end{align*}
\end{proof}

In \cite{alpay2020generalized, alpay2022generalized} the authors introduced a new class of functions called ML class which can be introduced as follows:
\begin{definition}[ML class]
We denote by ML the class of all entire functions $\varphi(z)=\displaystyle \sum_{n=0}^{\infty} a_nz^n$
satisfying the following conditions:
\begin{itemize}
\item[i)] $\varphi(0)=1$ and $\varphi'(0)>0;$
\item[ii)] $\varphi(z\overline{w})$ is a positive definite function on $\mathbb{C}$;
\item[iii)] $\varphi(-||.||^2/2)$ is a positive definite function on the space of Schwartz test functions $\mathcal{S}(\mathbb{R})$.
\end{itemize}
\end{definition}
\begin{remark}
If $\varphi$ is a function in the ML class then the map defined by $$T:s\in \mathcal{S}(\mathbb{R})\longrightarrow T(s):=\varphi(-||s||^2/2)$$ is continuous in the $L^2$ topology, and hence in the $\mathcal{S}$-topology, since the $L^2$ norm belongs to the set of norms defining the topology of the Fr\'echet space $\mathcal{S}$.
\end{remark}
\begin{theorem}
For a fixed $w\in \mathbb{C}$, the kernel function $\varphi_w=\eta_0\mathsf{K}_w$ defined on $\mathbb{C}$ by \eqref{Ker}, i.e: $$\varphi_w(z)=\eta_0 \mathsf{K}(z,w):=\eta_0\mathsf{E}(z\overline{w}), \quad \forall z\in \mathbb{C},$$ is positive definite.
\end{theorem}
\begin{proof}
We note that $\varphi_w(0)=\eta_0 \frac{1}{\eta_0}=1$; moreover we have
$$\displaystyle \varphi_w'(z)=\eta_0\sum_{k=1}^{\infty} \frac{k}{\eta_{k}}z^{k-1}\overline{w}^k,\quad \forall z,w\in\mathbb{C}.$$
Thus, we have also $\displaystyle \varphi_w'(0)=\frac{\eta_0}{\eta_1}>0.$
In order to justify that the kernel function $\mathsf{K}$ is positive definite, it is enough to observe that the function $\mathsf{E}(z\overline{w})$ can be factorized, so that for every $z,w\in \mathbb{C}$ we have
\begin{align*}
\displaystyle \mathsf{E}(z\overline{w})&=\int_\mathbb{R}A_z(x)\overline{A_w(x)}dx\\
&=\langle A_z, A_w \rangle_{L^2(\mathbb{R})}.\\
\end{align*}
\end{proof}
One of the most important results on ML functions is given by Theorem 2.5 in \cite{alpay2020generalized} since it allows to use a Bochner-Minlos type theorem leading to the following problem:
\begin{problem}
Can we prove that the function $\varphi_w(-||.||^2/2)$ associated to $\mathsf{E}$ is positive definite on $\mathcal{S}(\mathbb{R})$ and that therefore the function $\varphi_w$ belongs to the ML class ?
Is there a probability measure $P_\mathsf{E}$, such
that
\begin{equation}
\displaystyle \mathsf{E}\left(-\frac{||s||^2}{2}\right)=\int_{\mathcal{S}'(\mathbb{R})}e^{i\langle w,s \rangle}dP_{\mathsf{E}}(w); \quad s\in \mathcal{S}(\mathbb{R}) ?
\end{equation}
\end{problem}

An interesting example of ML functions which we propose here is inspired by the special function $E_n(z)$ given by formula \eqref{Enz}, which satisfies this special identity (see \cite[page 230]{abramowitz1964handbook})
\begin{equation}\label{ESF}
\displaystyle \int_0^\infty e^{-at}E_n(t)dt=\frac{(-1)^{n-1}}{a^n}\left(\log(1+a)+\sum_{k=1}^{n-1}\frac{(-a)^k}{k}\right); \quad a>-1.
\end{equation}
In order to present this example we denote by $\mathbb{D}=\lbrace{z\in \mathbb{C}, |z|<1}\rbrace$ the complex unit disk and consider on it the following special function:
\begin{definition}
For every $n=1,2,...$ we consider the function defined by
\begin{equation}
\phi_n(z)=\frac{1}{z^n}\sum_{k=n}^\infty \frac{z^k}{k}=\sum_{k=0}^\infty \frac{z^k}{k+n};\quad \forall z\in\mathbb{D}.
\end{equation}
\end{definition}

\begin{proposition}
For every $n=1,2,...$ we have
$$\phi_n(0)=\frac{1}{n}, \quad \phi'_n(0)=\frac{1}{n+1}.$$
Moreover, for every $p=0,1,...$ we have
$$\phi^{(p)}_n(0)=\frac{p!}{n+p}.$$
\end{proposition}
\begin{proof}
It is clear that the function $\phi_n(z)=\sum_{p=0}^{\infty} a_pz^p$ with coefficients $\displaystyle a_p=\frac{1}{p+n}$ is an analytic function. Moreover, using the Taylor series expansion we know that its coefficients are given by $$a_p=\frac{\phi_{n}^{(p)}(0)}{p!};\quad \forall p=0,1,\ldots$$, so that $\displaystyle \phi_{n}^{(p)}(0)=\frac{p!}{n+p}$ for every $p=0,1 \ldots$

\end{proof}
Now, for every $n=1, 2, \cdots,$ we define on the unit disk $\mathbb{D}$ the kernel function given by
\begin{equation}\label{HF}
k_n(z,w):=\phi_n(z\overline{w})=\sum_{k=0}^\infty \frac{z^k \overline{w}^k}{k+n}, \quad \forall z,w\in \mathbb{D}.
\end{equation}
It turns out that the kernel function in \eqref{HF} is related to the well-known Lerch transcendent which is a generalization of the Hurwitz zeta function $\zeta(s,a)$ defined by
\begin{equation}
\zeta(s,a):=\sum_{k=0}^\infty \frac{1}{(k+a)^s}, \quad Re(s)>1.
\end{equation} The Hurwitz function has the following important integral representation
\begin{equation}
\displaystyle \zeta(s,a)=\frac{1}{\Gamma(s)}\int_0^\infty \frac{t^{s-1}}{e^{at}(1-e^{-t})}dt; \quad Re(s)>1, Re(a)>0.
\end{equation}
We recall also a generalized Hurwitz function, sometimes called also the Lerch zeta transcendent function or (Hurwitz–Lerch zeta function)  defined by
\begin{equation}
\Phi(z,s, a ):=\sum_{k=0}^\infty \frac{z^k}{(k+a)^s}, \quad |z|<1
\end{equation}
and which has the following integral representation
\begin{equation}
\displaystyle \Phi(z,s, a )=\frac{1}{\Gamma(s)}\int_0^\infty \frac{t^{s-1}e^{-at}}{1-ze^{-t}}dt.
\end{equation}
So, in our particular case related to the kernels $k_n$ we have
\begin{equation}
k_n(z,w)=\phi_n(z\overline{w})=\sum_{k=0}^\infty \frac{z^k \overline{w}^k}{k+n}=\Phi(z\overline{w},1,n).
\end{equation}
It is important to observe that the kernel function $k_1$  on the unit disk $\mathbb{D}$ corresponds to the kernel of the classical Dirichlet space $\mathcal{D}$ (see Theorem 1.2.3 of the book \cite{el2014primer} ). This follows from the fact that \begin{equation}
k_1(z,w)=\phi_1(z\overline{w})=\sum_{k=0}^\infty \frac{z^k\overline{w}^k}{k+1}=\frac{1}{z\overline{w}}\log\left(\frac{1}{1-z\overline{w}}\right), \quad z, w\in \mathbb{D}\setminus \lbrace 0 \rbrace.
\end{equation}
\begin{proposition}
The function $\phi_n(z\overline{w})$ is positive definite on $\mathbb{D}$.
\end{proposition}
\begin{proof}
We note that the function $\phi_1(z\overline{w})$ is positive definite since it is the reproducing kernel of the classical Dirichlet space $\mathcal{D}$. Moreover, setting
$$g_{k,n}(z):=\frac{z^k}{\sqrt{k+n}}, \quad k=0,1,\cdots, n=1, 2\cdots $$
we obtain
$$k_n(z,w)=\phi_n(z\overline{w})=\sum_{k=0}^\infty g_{k,n}(z)\overline{g_{k,n}(w)}.$$
So the kernel $k_n(z,w)=\phi_n(z\overline{w})$ is positive definite.
\end{proof}
\begin{proposition}
For every $a>-1$ we have

\begin{equation}\label{F1}
\phi_n(-a)=\displaystyle \int_0^\infty e^{-at}E_n(t)dt.
\end{equation}

Moreover, for every $s\in\mathcal{S}(\mathbb{R})$ we have
\begin{equation}\label{F2}
\phi_n\left (-\frac{||s||^2}{2}\right)=\displaystyle \int_0^\infty e^{-\frac{||s||^2}{2} t}E_n(t)dt.
\end{equation}

\end{proposition}
\begin{proof}
Let $a>-1$, we first recall that

$$\log(1+a)= -\sum_{k=1}^{\infty}\frac{(-1)^ka^k}{k}.$$

Then, by definition of $\phi_n$ we have
\begin{align*}
\displaystyle \phi_n(-a)&=\frac{(-1)^n}{a^n} \sum_{k=n}^{\infty}\frac{(-a)^k}{k}\\
&=\frac{(-1)^n}{a^n}\left( \sum_{k=1}^{n-1}\frac{(-a)^k}{k}+  \sum_{k=n}^{\infty}\frac{(-a)^k}{k}-  \sum_{k=1}^{n-1}\frac{(-a)^k}{k}\right)\\
&=\frac{(-1)^n}{a^n}\left( \sum_{k=1}^{\infty}\frac{(-a)^k}{k}-  \sum_{k=1}^{n-1}\frac{(-a)^k}{k}\right)\\
&= \frac{(-1)^{n-1}}{a^n}\left( \log(1+a)+  \sum_{k=1}^{n-1}\frac{(-a)^k}{k}\right)\\
\end{align*}
We then use formula \eqref{ESF} to deduce that

$$\phi_n(-a)=\displaystyle \int_0^\infty e^{-at}E_n(t)dt.$$

Now, let $s\in\mathcal{S}(\mathbb{R})$, in order to justify the formula \eqref{F2} we just need to set $a=||s||^2/2$ and insert it in formula \eqref{F1} to conclude the proof.

\end{proof}

\begin{theorem}
Consider the function defined by $\tilde{\phi_n}(z)=n\phi_n(z)$ for every $z\in \mathbb{C}$. There exists a uniquely defined probability measure $P_{\phi_n}$, such
that
\begin{equation}
\displaystyle\tilde{ \phi_n}\left(-\frac{||s||^2}{2}\right)=\int_{\mathcal{S}'(\mathbb{R})}e^{i\langle w,s \rangle}dP_{\tilde{\phi_n}}(w); \quad s\in \mathcal{S}(\mathbb{R}).
\end{equation}
\end{theorem}
\begin{proof}
We have $\tilde{\phi_n}(0)=n\phi_n(0)=1,$
and we can see that the function $\tilde{\phi_n}$ belongs to the class ML. So,
applying Theorem 2.5 of \cite{alpay2020generalized} for ML functions we have the result.
\end{proof}
\section*{Acknowledgments}
Daniel Alpay thanks the Foster G. and Mary McGaw Professorship in Mathematical Sciences, which supported this research. Kamal Diki thanks the Grand Challenges Initiative (GCI) at Chapman University, which supported this research. 

\end{document}